\documentclass[reqno]{amsart}
\usepackage{amsmath}

\newtheorem{theorem}{Theorem}[section]
\newtheorem{lemma}[theorem]{Lemma}
\newtheorem{proposition}[theorem]{Proposition}

\theoremstyle{definition}
\newtheorem{definition}[theorem]{Definition}

\theoremstyle{remark}

\numberwithin{equation}{section}

\newcommand{\C}{{\mathbb C}}
\newcommand{\be}[1]{\begin{equation}\label{#1}}
\newcommand{\ee}{\end{equation}}

\begin{document}

\title[Discrete Connection Laplacians]{Discrete Connection Laplacians}
\author{SVETOSLAV ZAHARIEV}
\address {Department of Mathematics and Computer Science, Lehman College of CUNY,
 250 Bedford Park Boulevard West, Bronx, NY 10468,
U.S.A.} \email {szahariev@gc.cuny.edu}


\subjclass[2000]{Primary 58J50}

\date{}

\dedicatory{}

\commby{Mikhail Shubin}

\begin{abstract}
To every Hermitian vector bundle with connection over a compact
Riemannian manifold $M$ one can associate a corresponding connection
Laplacian acting on the sections of the bundle. We define analogous
combinatorial, metric dependent Laplacians associated to
triangulations of $M$ and prove that their spectra converge, as the
mesh of the triangulations approaches zero, to the spectrum of the
connection Laplacian.
\end{abstract}

\maketitle

\section{Introduction}

Let $M$ denote a compact Riemannian manifold. Combinatorial
analogues of the Laplacian on $M$ have been considered in various
contexts since the appearance of \cite{E}.  In the 1970s Dodziuk
defined combinatorial, metric dependent Laplacians associated to
triangulations of $M$ and proved that under certain technical
assumptions the spectra of these operators converge to the spectrum
of the Laplacian acting on functions on $M$ as the mesh of the
triangulations tend to zero, cf. \cite{D}. This result was extended
to the Laplacian acting on differential forms in \cite{DP} and used
in one of proofs of the Ray-Singer conjecture, cf. \cite{Mu}.

A natural generalization of the Laplacian on $M$ is the connection
(also called rough) Laplacian in whose definition the exterior $d$
is replaced by the  covariant derivative $\nabla$ associated to a
connection on a Hermitian vector bundle over $M$. Not much seems to
be known about the spectra of such operators in general (but see
\cite{BBC}). When the bundle is trivial of rank 1 (and usually in a
noncompact setting) such connection Laplacians are also known and
studied in the mathematical physics literature, e.g. \cite{MSh}, as
magnetic Schr\"{o}dinger operators.

In this note we propose two definitions of combinatorial analogs of
such connection Laplacians. Both rely on replacing the wedge product
of differential forms by a graded commutative but nonassociative cup
product on simplicial cochains, which is probably part of folklore,
and whose definition can be found for example in \cite{BR}. In the
first discretization scheme, which is applicable only for trivial
bundles, we define combinatorial Laplacians corresponding to a
coboundary operator twisted by a 1-cochain  defined using the above
cup product. The convergence of eigenvalues in degree zero (i.e. for
the Laplacian acting on sections only rather than on vector bundle
valued forms) in this case is the content of Theorem \ref{conver1}.
The second construction of discrete connection Laplacian uses the
well known fact, cf. \cite{NR}, that any bundle with connection can
be embedded into a trivial bundle in such a way that the connection
is induced by the standard connection on the trivial bundle. Then
the idea is to discretize this "classifying map". We establish the
corresponding convergence of spectra in Theorem \ref{convgeneral}.

We mention that our approach leads to a natural generalization of
the discrete magnetic Laplacians on graphs studied in \cite{S} and
\cite{MY} to simplicial complexes. One can expect that an analogous
convergence of spectrum distribution functions result for such
twisted Laplacians can be proved also for noncompact universal
covering manifolds, using the method of \cite{CCMP}. Finally, we
note that a quite different discretization scheme for the degree
zero case can be found in \cite{M1} and \cite{M2}.

\section{The de Rham and Whitney maps}

In this note $M$ will always denote  a  closed oriented Riemannian
manifold of dimension $N$. Let $\Omega(M)$ be the vector space of
all smooth complex valued differential forms on $M$ and
$L^{2}(\Lambda(M))$, the space of square integrable forms on $M$
with respect to the volume element induced by the metric. Let $K$ be
a smooth triangulation of $M$ and $C(K)$ be the vector space of
(oriented) cochains of $K$ with complex coefficients. In what
follows, $\|\omega\|$ will always denote the $L^{2}$-norm of the
form $\omega$ and $\|c\|$, the norm of the cochain $c$ given by the
canonical inner product on cochains with respect to which simplices
are orthogonal.

Recall that for every $q>0$ the corresponding de Rham map $R^{K}:
\Omega^{q}(M)\rightarrow C^{q}(K)$ is given by $
R^{K}(\omega)(\sigma)=\int_{\sigma}\omega$, where $\omega \in
\Omega^{q}(M)$ and $\sigma$ is a $q$-simplex. When $q=0$ the de Rham
map is defined to be simply the evaluation of a function on the
$0$-skeleton of $K$.

The Whitney mapping $W^{K}$ goes in the reverse direction and is
defined as follows. Let $\mu_{i}=\mu_{v_{i}}$ denote the barycentric
coordinate function corresponding to a vertex $v_{i}$ and let
$\sigma=[v_{0},....v_{q}]$ be a $q$-simplex, where $q>0$. We set

\[W^{K}(\sigma)=q!\sum_{i=0}^{q}(-1)^{i}\mu_{i}d\mu_{0}\wedge
...\wedge d\mu_{i-1}\wedge d\mu_{i+1}\wedge ...\wedge d \mu_{q}.\]

If $q=0$ we set $W^{K}(v_{0})=\mu_{0}$. This expression defines a
map: $C^{q}(K)\rightarrow L^{2}(\Lambda^{q}(M))$. We note that the
image of the Whitney map does not consist of smooth forms and
summarize its basic properties in the following

\begin{proposition}\label{Wprop}
i) On the complement of the $(N-1)$-skeleton of $K$ one has
$W^{K}d^{K}=dW^{K}$, where $d^{K}$ denotes the simplicial coboundary
operator of $K$.

(ii)Let $j_{\sigma}$ denote the inclusion map: $\sigma
\hookrightarrow M$ and $j^{*}_{\sigma}$ the corresponding pull-back
map on forms. Then the values of $j^{*}_{\sigma}W^{K}c$ (for $c$ an
arbitrary cochain) depend  only on the value of $c$ at $\sigma$.
Therefore $R^{K}$ is well-defined on the image of $W^{K}$ and one
has $R^{K}W^{K}=id$ on $C(K)$.

(iii) The support of $W^{K}(\sigma)$ is contained in the closed star
of $\sigma$.
\end{proposition}

 For  proofs the reader is referred to \cite{Wh} or
\cite{D}.

The injectivity of the Whitney map allows us to define a new, metric
dependent inner product on the space of cochains, namely, for
$c_{1},c_{2} \in C(K)$ we set $\langle c_{1},c_{2}\rangle
_{W}=\langle W^{K}c_{1},W^{K}c_{2}\rangle_{L^{2}}$. We will refer to
this inner product as the {\em Whitney inner product}.

\section{The trivial bundle case}
\subsection{A commutative cup product}
We first let $K$ be any finite simplicial complex.
\begin{definition}[ \cite{BR}, \cite{A}, \cite{W}]\label{cupdef}
 For two oriented simplices
$\sigma_{1} \in C^{p}(K)$ and $\sigma_{2} \in C^{q}(K)$, set
$\sigma_{1} \cup \sigma_{2}=0$ unless $\sigma_{1}$ and $\sigma_{2}$
meet at precisely one vertex and span a ($p+q$)-simplex $\tau$, in
which case define $ \sigma_{1} \cup \sigma
_{2}=\epsilon(\sigma_{1},\sigma_{2} )\frac{p!q!}{(p+q+1)!}\tau$,
where the sign $\epsilon(\sigma_{1},\sigma_{2} )=\pm 1$ is
determined by the equation $orientation (\sigma_{1}).orientation
(\sigma_{2})=\epsilon(\sigma_{1},\sigma_{2}) .orientation (\tau)$.
\end{definition}

One checks that this defines a graded commutative non-associative
bilinear operation on $C(K)$ with respect to which the coboundary
operator is a derivation.

In the case when $K$ is the underlying complex of a smooth
triangulation of a compact manifold $M$ it was observed in \cite{Du}
( see also \cite{A} and \cite{W} for more details) that this cup
product admits an alternative description in terms of the Whitney
and De Rham maps. Namely, for $a,b \in C(K)$ one has

\begin{equation}\label{cup}
a \cup b=R^{K}(W^{K}a \wedge W^{K}b)
\end{equation}

Now let  $A$ be a fixed real smooth 1-form. Consider the {\em
twisted} exterior differential $d_{A}=d+iA\wedge .$ and the
corresponding Laplacian
$\triangle_{A}=d^{*}_{A}d_{A}+d_{A}d^{*}_{A}$. We use the above to
define discrete analogues of these operators.

Let $a$ be a fixed 1-cochain. Define the {\em twisted} coboundary
operator associated to it as $d^{K}_{a}=d^{K}+ia\cup .$, where
$d^{K}$ is the usual coboundary of the simplicial complex $K$.
Define the {\em twisted discrete Laplacian} by
$\triangle^{K}_{a}=\left(d^{K}_{a}\right)^{*}d^{K}_{a}+d^{K}_{a}\left(d^{K}_{a}\right)^{*}$,
where the adjoint is taken with respect to the inner product given
by the Whitney map.

From now on we consider only triangulations $K$ which are
subdivisions of a fixed triangulation and whose fullness is bounded
away from $0$ (see \cite{Wh} or \cite {DP} for a definition). Then
the cup product on cochains introduced above approximates the wedge
product on forms according to the following result established in
\cite{W}.(The case when $\omega_{2}=1$ was proved in \cite{D}.)

\begin{theorem}\label{cupth}
Let $\sigma$ be a simplex in $K$, $p$ - any point in the interior of
$\sigma$ and $x_{1},...,x_{N}$ - local coordinates around $\sigma$.
Let $\omega_{1},\omega_{2}\in \Omega(M)$. There exists constant
$C>0$ independent of $\omega_{1},\omega_{2},\sigma$ and $K$ such
that

\begin{multline}\label{cupest}
 \left|W^{K}\left(R^{K}\omega_{1}\cup
R^{K}\omega_{2}\right)(p)-\omega_{1}\wedge
\omega_{2}(p)\right|_{p}\\
\leq
Ch\left(\sup|\omega_{1}|\sup\left|\frac{\partial\omega_{2}}{\partial
x_{i}}\right|+\sup|\omega_{2}|\sup\left|\frac{\partial\omega_{1}}{\partial
x_{i}}\right|\right).
 \end{multline}

Here the suprema are taken over all partial derivatives and all
points in $\sigma$.

\end{theorem}
Let us recall the definition of Sobolev spaces of vector bundle
valued forms. Let $E$ be a Hermitian complex vector bundle over $M$
and let $\Omega(M,E)$ denote the smooth differential forms on $M$
with values in $E$. We fix a Hermitian connection on $E$ which
defines a covariant differential $\nabla:\Omega ^{*}(M,E)\rightarrow
\Omega ^{*+1}(M,E)$.
 Let
$H^{r}(\Lambda(M,E))$ be the completion of $\Omega(M,E)$ with
respect to the norm
$\|\omega\|_{r}=\int_{M}|(1+\triangle_{E})^{\frac{r}{2}}\omega
|^{2}dvol$, where $\triangle _{E}=\nabla ^{*}\nabla + \nabla \nabla
^{*}$, $\nabla ^{*}$ being the formal adjoint of $\nabla$, is the
connection Laplacian associated to $E$ and $\nabla$.

When $E$ is the trivial line bundle we write $H^{r}(\Lambda(M))$ for
$H^{r}(\Lambda(M,E))$. In this case, integrating the estimate in
Theorem \ref{cupth} with $\omega_{2}=1$ and using the standard
Sobolev inequality (see e.g. \cite{N}) one obtains for every
$r>\frac{N}{2}+1$

\begin{equation}\label{best}\left\|W^{K}R^{K}\omega - \omega \right\|\leq Ch
\|\omega\|_{r}.\end{equation}

\subsection{Convergence of eigenvalues}

We set $a=R^{K}A$ and use the results above to compare the operators
$d_{A}$ and $d^{K}_{a}$. It turns out that the analogue of the
identity from Proposition \ref{Wprop} (i) holds  approximately in
the twisted setting.

\begin{lemma}\label{lemmamainest1}
For every $\omega \in \Omega (M)$ and every point $p \in M$ we have

\begin{equation}\label{est1} \left| W^{K}d^{K}_{a}R^{K}\omega - d_{A}W^{K}R^{K}\omega
\right|_{p}\leq C_{A,\omega}h,
 \end{equation}

\begin{equation}\label{est2} \left| W^{K}d^{K}_{a}R^{K}\omega - W^{K}R^{K}d_{A}\omega
\right|_{p}\leq C_{A,\omega}h,
\end{equation}

where $C_{A,\omega}$ is a constant depending only on $A$ and
$\omega$ and their first derivatives.

\end{lemma}

\begin{proof}

In order to show (\ref{est1}), we write
\begin{multline*}W^{K}d^{K}_{a}R^{K}\omega - d_{A}W^{K}R^{K}\omega =
W^{K}\left(R^{K}iA \cup R^{K}\omega \right)- iA\wedge
W^{K}R^{K}\omega \\=W^{K}\left(R^{K}iA \cup R^{K}\omega \right)-iA
\wedge \omega - iA \wedge \left(W^{K}R^{K}\omega -\omega \right)
\end{multline*}
and then apply  Theorem \ref{cupth} and the basic estimate
(\ref{best}). Similarly, observe that
\begin{multline*}
W^{K}d^{K}_{a}R^{K}\omega - W^{K}R^{K}d_{A}\omega =
W^{K}\left(R^{K}iA \cup R^{K}\omega \right) - W^{K}R^{K}(iA \wedge
\omega)\\
=W^{K}\left(R^{K}iA \cup R^{K}\omega \right)-iA \wedge \omega + iA
\wedge \omega - W^{K}R^{K}(iA \wedge \omega).
\end{multline*}
This easily implies (\ref{est2}).
\end{proof}

Note that $d_{A}$ extends to a map: $H^{1}(\Lambda(M))\rightarrow
L^{2}(\Lambda(M))$, denoted also by $d_{A}$.

\begin{proposition}\label{mainest1}
i) For every $\omega \in \Omega (M)$  we have

\be{est3} \left\|W^{K}d^{K}_{a}R^{K}\omega-d_{A}\omega \right\|\leq
C_{A,\omega}h,\ee

where $C_{A,\omega}$ is a constant depending only on $A$ and
$\omega$ and their derivatives.

ii) For every $c \in C^{0}(K)$ we have

\be{est4} \left\|W^{K}d^{K}_{a}c-d_{A}W^{K}c \right\| \leq
C_{A}\left\|W^{K}c \right\|_{1}h .\ee

As above,  $C_{A}$  depends only on $A$  and its derivatives of
order up to $r$.

\end{proposition}

\begin{proof}

 The estimate (\ref{est3}) is an easy
consequence of (\ref{est2}). Observe that the image of $C^{0}(K)$
under the Whitney map consists of continuous piecewise linear
functions and therefore is in $H^{1}\left(\Lambda(M)\right)$.
Denoting by $\left\|.\right\|_{r}^{\sigma}$ for the $r$-th Sobolev
norm on a closed $N$-simplex $\sigma$ computed using local
coordinates in a neighborhood of $\sigma$, we obtain

\begin{multline*}
\left\|W^{K}d^{K}_{a}c-d_{A}W^{K}c \right\|^{2} \leq \sum_{\sigma
\in
K}\int_{\sigma}\left|W^{K}d^{K}_{a}c-d_{A}W^{K}c\right|^{2}dvol\leq
\\
\sum_{\sigma \in
K}C'_{A,h}h^{N+2}\left(\|W^{K}c\|_{r}^{\sigma}\right)^{2}\leq
\sum_{\sigma \in
K}C_{A}h^{2}\left(\|W^{K}c\|_{1}^{\sigma}\right)^{2}\leq
C_{A}h^{2}\|W^{K}c\|_{1}^{2}.
\end{multline*}

Above we used the local estimate (\ref{est1}) with $\omega=W^{K}c$,
a version of the Sobolev inequality which can be applied to each
closed $N$-simplex (see, e.g. \cite[Theorem 3.9]{Ag}) and the
inequality $C'_{A,h}\leq C_{A}h^{N}$.

\end{proof}

It is well-known that $\triangle_{A}$ with domain all smooth forms
is an essentially selfadjoint elliptic operator on the compact
manifold $M$ and hence its $L^{2}$ closure, denoted also by
$\triangle_{A}$, has purely discrete spectrum. Let $\lambda_{0}\leq
\lambda_{1} \leq \lambda_{2} \leq \ldots $ be the eigenvalues,
repeated according to their multiplicities, of $\triangle _{A}$
acting on functions and let $\lambda ^{K}_{0}\leq \lambda ^{K}_{1}
\leq \ldots \leq \lambda ^{K}_{\dim C^{0}(K)}$ denote the
eigenvalues of $\triangle ^{K}_{a}$ acting on 0-cochains. We have
the following convergence of eigenvalues result.

\begin{theorem}\label {conver1}
Let $j \leq \dim C^{0}(K)$. There exist positive constants $C_{A}$
and $C_{j,A}$ such that following inequalities hold.

\[\lambda ^{K}_{j}-C_{j,A}h\leq \lambda_{j}\leq \frac{\left(\sqrt{\lambda^{K}_{j}}+
C_{A}h\right)^{2}}{1-C_{A}h}.\]

\end{theorem}

\begin{proof}
We begin by proving the second inequality. It follows from
(\ref{est4}) that for every 0-cochain $c$ and some constant
$C_{A}\geq 0$ depending only on $A$ and its first derivatives one
has

\[ \left\|W^{K}d^{K}_{a}c-d_{A}W^{K}c \right\| \leq
C_{A}h\left(\left\|W^{K}c
\right\|+\left\|d_{A}W^{K}c\right\|\right).\]

Thus we obtain

 \be{theoest1} \left\|W^{K}d^{K}_{a}c\right\|\geq
\left(1-C_{A}h\right)\left\|d_{A}W^{K}c\right\|-C_{A}h\left\|W^{K}c\right\|.\ee

Using this, the min-max principle (see e.g. \cite{RS}) gives

\begin{multline*}
\sqrt{\lambda^{K}_{j}}=\sup_{c_{1}\ldots,c_{j-1}\in
C^{0}(K)}\inf_{\substack{c\neq 0, \langle c,c_{k}\rangle =0,\\
k=1,2,\ldots,j-1}}
\frac{\left\|W^{K}d^{K}_{a}c\right\|}{\left\|W^{K}c \right\|}\\
\geq \sup_{c_{1},\ldots,c_{j-1}\in C^{0}(K)}\inf_{\substack{ c\neq
0, \langle c,c_{k}\rangle =0,
\\ k=1,2,\ldots,j-1}}\left(\frac{\left(1-C_{A}h
\right)\left\|d_{A}W^{K}c\right\|}{\left\|W^{K}c
\right\|}-C_{A}h\right)\\
\geq \left(1-C_{A}h \right)\sup_{f_{1},\ldots,f_{j-1}\in
W\left[C^{0}(K)\right]}\inf_{\substack{f \in
W\left[C^{0}(K)\right]\setminus \{0\},\\ \langle f,f_{k}\rangle
=0,k=1,...,j-1}}\frac{\left\|d_{A}f\right\|}{\left\|f\right\|}-C_{A}h\\
\geq\left(1-C_{A}h\right)\sqrt{\lambda_{j}}-C_{A}h.
\end{multline*}

Hence the second inequality holds.

To establish the first inequality, denote by $V_{j}$  the vector
space spanned by the first $j$ eigenfunctions of $\triangle_{A}$.
Then it follows from (\ref{est3}) that for each $f \in
V_{j}\setminus \{0\}$ and for some constant $C_{j,A}$ depending only
on $j$ and $A$ we have

\[\left|\frac{\langle W^{K}d^{K}_{a}R^{K}f,W^{K}d^{K}_{a}f \rangle}
{\langle W^{K}R^{K}f,W^{K}R^{K}f \rangle}-\frac{\langle
d_{A}f,d_{A}f\rangle}{\langle f,f \rangle}\right|\leq C_{j,A}h.\]

Then the proof proceeds in exactly the same fashion as in
\cite[Theorem 5.7]{D}.

\end{proof}

\section{The general bundle case}

\subsection{The universal connection}
 \indent Now let $E$ be a Hermitian complex vector bundle
of rank $d$ over $M$. We denote  the space of $L^{2}$ forms with
values in $E$ defined using the inner product on $E$ by
$L^{^{2}}\left(\Lambda(M),E \right)$ and fix a Hermitian connection
on $E$ which defines a covariant differential $\nabla:\Omega
^{*}(M,E)\rightarrow \Omega ^{*+1}(M,E)$. We shall need the
following theorem on the existence of universal connections due to
Narasimhan and Ramanan \cite{NR}, in the form proved in the appendix
of \cite{Q}.

\begin{theorem}\label{univ}

There exist a trivial Hermitian bundle $V$ and an isometric
embedding of bundles $i: E \hookrightarrow V$ such that
$\nabla=i^{*}\circ d \circ i$, where $i^{*}$ denotes the fiberwise
adjoint taken with respect to the corresponding inner products.

\end{theorem}

Let us recall briefly how the Narasimhan-Ramanan theorem is proved
as the construction will be needed in the sequel. One first shows
that the statement is true for trivial bundles. Then one chooses a
finite open cover $\{U_{l}\}$ of $M$ such that $E$ restricted to
every $U_{l}$ is trivial and a partition of unity $\{\psi_{l}\}$
subordinate to this cover such that $\sum_{l}\psi_{l}^{2}=1$. If
$i_{l}$ denotes the already constructed isometric embedding from
$E|_{U_{l}}$ to a trivial bundle $V_{l}$ one defines
\begin{equation}\label{defofi} i=\sum_{l}i_{l}\psi_{l}:E \rightarrow V,
\end{equation}
where $V=\bigoplus _{l}V_{l}$.

  In what follows, we will assume that
we have fixed a triangulation $K_{0}$ fine enough so that there
exists a finite open cover $\{U_{l}\}$  as above with the additional
properties that (1) each  $\overline{U}_{l}$ is a subcomplex of
$K_{0}$, and (2) $E$ restricted to  $\overline{U}_{l}$ is trivial,
and we will consider only subdivisions $K$ of $K_{0}$. We will also
assume that the map $i$ is constructed using such an open cover.
\subsection{Twisted de Rham and Whitney maps}
We first define the space of cochains on which the discrete
connection Laplacian will act.
 We denote the induced map from
 $L^{^{2}}\left(\Lambda(M),E \right)$ to $L^{^{2}}\left(\Lambda(M),V
\right)$ also by $i$. We would like to define cochains on $K$ with
values in the bundle $E$. To this end, fix a reference point
$p_{\sigma}$ in the interior of each simplex $\sigma$ in $K$.

\begin{definition}\label{deftwist}
 Let $C(K,E)$, {\em the twisted cochains with values in $E$}, denote the set of all maps from
the set of all simplices in $K$ to the total space of $E$ which
assign to each simplex $\sigma$ a vector in the fibre
$E_{p_{\sigma}}$.
\end{definition}

 This is a vector space  with an inner product induced
by the Hermitian structure on $E$. The restriction of $i$ to fibers
defines a map, denoted by $i^{K}$, from $C(K,E)$ to $C(K,V)$ which
we identify with the cochains taking values in $\C ^{n}$, where $n$
is the rank of the trivial bundle $V$.

We shall also need the following notation. Let $I_{l}^{K}:
C(\overline{U_{l}},\C^{d})\rightarrow C(\overline{U_{l}},\C^{n})$ be
the operator given by $(I_{l}^{K}c)_{p}=\sum_{s}R^{K}(i_{ps})\cup
c_{s}$. Here $i_{ps}$ are the entries of the matrix of $i$ and
$c_{p}$ the components of the vector valued cochain $c$. We define
the operator $\Psi_{l}^{K}$ by
$\Psi_{l}^{K}c_{p}=R^{K}(\psi_{l})\cup c_{p}$. Finally we define
$I^{K}: C(K,E)\rightarrow C(K,V)$ by
$I^{K}c=\sum_{l}I_{l}^{K}\Psi_{l}^{K}c$. Now define an operator
acting on $C(K,E)$, the discretization of $\nabla$, by setting
$\nabla ^{K}=\left(I^{K}\right)^{*}d^{K}I^{K}$.

In order to compare $\nabla ^{K}$ and $\nabla$ we have to introduce
the appropriate {\em twisted Whitney and de Rham maps}. We set
$\widetilde{W}^{K}=i^{*}W^{K}I^{K}$ and
$\widetilde{R}^{K}=\left(I^{K}\right)^{*}R^{K}i$. Here $W^{K}$ and
$R^{K}$ are the usual Whitney and de Rham maps but now acting
between the spaces of vector valued cochains and differential forms,
$C(K,V)$ and $L^{^{2}}\left(\Lambda(M),V \right)$, respectively.

It is easily seen that $\widetilde{W}^{K}$ is injective on
0-cochains. More generally, we have:

\begin{lemma}\label{inject}
For sufficiently fine triangulations $K$ the map $\widetilde{W}^{K}$
is injective on $C^{q}(K,E)$ for all $q$.
\end{lemma}

The proof will be given in the next section.  We can now define the
corresponding Whitney inner product on $C(K,E)$ as in Section 2 and
the {\em discrete connection Laplacian} $\triangle
^{K}_{E}=\left(\nabla^{K}\right)^{*}\nabla ^{K}+ \nabla
^{K}\left(\nabla ^{K}\right)^{*}$, where the adjoint is taken with
respect to this inner product. Let  $\lambda_{0}\leq \lambda_{1}
\leq \lambda_{2} \leq \ldots $ denote the eigenvalues repeated
according to their multiplicities, of  (the $L^{2}$ closure of) the
connection Laplacian $\triangle _{E}=\nabla ^{*}\nabla + \nabla
\nabla ^{*}$ acting on sections of $E$, and let $\lambda
^{K}_{0}\leq \lambda ^{K}_{1} \leq \ldots \leq \lambda ^{K}_{\dim
C^{0}(K,E)}$ denote the eigenvalues of $\triangle ^{K}_{E}$ acting
on 0-cochains with values in $E$. We have the following convergence
of spectra result.

\begin{theorem}\label{convgeneral}
There exist positive constants $C^{i}$ and $C^{i}_{j}$ such that
following inequalities hold for  $j \leq \dim C^{0}(K,E)$.
\[\lambda ^{K}_{j}-C^{i}_{j}h\leq \lambda_{j}\leq \frac{\left(\sqrt{\lambda^{K}_{j}}+
C^{i}h\right)^{2}}{1-C^{i}h}.\]
\end{theorem}

\section{Proof of  Theorem \ref{convgeneral}}\label{sec6}
\subsection{Estimate from below}

In this subsection we estimate from below the eigenvalues of the
connection Laplacian in terms of the eigenvalues of its
discretizations, i.e., we prove the first inequality in Theorem
\ref{convgeneral}. We denote the projection $ii^{*}$ by $P$ and
introduce the "almost" projection operators
$P^{K}=I^{K}\left(I^{K}\right)^{*}$. Our next lemma states that the
Whitney and de Rham maps approximately interwine $P$ and $P^{K}$.

\begin{lemma}\label{lebelow}
 For $r>\frac{N}{2}+2$,
every form $\omega \in \Omega(M,V)$ and every vector valued function
$f \in \Omega ^{0}(M,V)$ one has
\begin{equation} \label{projR}
\left\|W^{K}P^{K}R^{K}\omega - W^{K}R^{K}P\omega \right\|\leq
C_{i}h\|w\|_{r},
\end{equation}
\begin{equation} \label{projW}
\left\|W^{K}P^{K}R^{K}\omega - PW^{K}R^{K}\omega \right\|\leq
C_{i}h\|w\|_{r},
\end{equation}
\begin{equation} \label{projRf}
P^{K}R^{K}f=R^{K}Pf.
\end{equation}
Here $C_{i}$ denotes a constant depending only on the map $i$ and
its derivatives.
\end{lemma}

\begin{proof}
We  first show that (\ref{projR}) holds.  In the course of the
proof, various constants depending on $i$ will all be denoted by
$C_{i}$.

We first assume that the bundle $E$ is trivial and identify $C(K,E)$
with $C(K,\C ^{d}$) and $\Omega (M,E)$ with $\Omega (M,\C ^{d})$.
Let $\sigma$ be a $q$-simplex in $K$ and let $\omega \in \Omega ^{q}
(M,E)$. We write $R^{K}$ (respectively $W^{K}$) for both de Rham
(Whitney) maps defined on $C(K,\C ^{d}$) and $C(K,V)$ (respectively
on $\Omega (M,\C ^{d})$ and $\Omega (M,V)$) and observe that in view
of \cite[Lemma 7.22]{DP} for every $c \in C(K,V)$ one has
\begin{equation}\label{DPlemma}
C_{1}h^{N-2q}\|c\|^{2}\leq\|W^{K}c\|^{2}\leq
C_{2}h^{N-2q}\|c\|^{2},\end{equation} where the constants $C_{1}$
and $C_{2}$ do not depend on  $K$. From the definitions we have $
I^{K}R^{K}\omega  - R^{K}i\,\omega = \sum_{s}\left(
R^{K}(i_{ps})\cup R^{K}\omega _{s}-R^{K}(i_{ps} \omega
_{s})\right)$. Thus, combining (\ref{best}), (\ref{cupest}) and
(\ref{DPlemma}) we conclude that
\begin{multline*}
 \left\|W^{K} I^{K}R^{K}\omega  - W^{K}R^{K}i\,\omega \right\|\\ \leq
\sum_{s}\left( \left\|W^{K}\left(R^{K}(i_{ps})\cup R^{K}\omega
_{s}-i_{ps}\omega_{s}\right)\right\|+\left\|W^{K}R^{K}i_{ps}\omega_{s}-i_{ps}\omega_{s}
\right\|\right).
\end{multline*}
Now  (\ref{DPlemma})implies that
\begin{equation}\label{Iiest}\left\| I^{K}R^{K}\omega  - R^{K}i\,\omega \right\|^{2} \leq
C_{i}h^{2q-N+2}\|\omega\|_{r}^{2}.\end{equation} It is easy to check
that $((I^{K})^{*}c)_{p}=\sum_{s}R^{K}(i^{*}_{ps})c_{s}$. Hence we
also have
\begin{equation}\label{Istarest}\left\| (I^{K})^{*}R^{K}\omega  - R^{K}i^{*}\omega \right\|^{2}
\leq C_{i}h^{2q-N+2}\|\omega\|_{r}^{2}.\end{equation} Next, observe
that the operators $I^{K}$ are uniformly bounded with respect to
 $K$.
 Indeed, from the combinatorial definition of the cup product and the
mean value theorem one easily obtains \be{Iandi}\left\|
i^{K}R^{K}\omega  - I^{K}R^{K}\omega \right\|^{2} \leq
C_{i}h^{2q-N+2}\|\omega\|_{r}^{2}\ee and $i^{K}$ is certainly
uniformly bounded. We use this fact in the following estimate.
\begin{multline*}\left\| P^{K}R^{K}\omega  - R^{K}P\omega
\right\|^{2} \\
\leq \left\|I^{K}(I^{K})^{*} R^{K}\omega -I^{K}R^{K}i^{*}\omega
\right\|^{2}+\left\|I^{K}R^{K}i^{*}\omega
-R^{K}ii^{*}\omega\right\|^{2}\leq
C_{i}h^{2q-N+2}\|\omega\|_{r}^{2}.
\end{multline*}
 Applying (\ref{DPlemma}) again, we see that
\[\left\| W ^{K}P^{K}R^{K}\omega  - W^{K}R^{K}P\omega \right\|
\leq C_{i}h\|\omega\|_{r}.\] Now let $E$ be an arbitrary bundle. We
use the open cover $\{U_{l}\}$ and the partition of unity
$\{\psi_{l}\}$ with the properties described after Theorem
\ref{univ} and define $P_{l}=i_{l}i_{l}^{*}$. Then, using
(\ref{defofi}), we find that
$P=\sum_{l,k}i_{l}\psi_{l}\psi_{k}i_{k}^{*}=
               \sum_{l}i_{l}\psi_{l}^{2}i_{l}^{*}$. On the other hand, $P^{K}$ is not equal to
$\sum_{l}I_{l}^{K}(\Psi_{l}^{K})^{2}(I _{l}^{K})^{*}$ due to the
non-associativity of our cup product. However, according to Theorem
\ref{cupth}, this product is approximately associative, therefore we
can replace $P^{K}$ by
 $\sum_{l}I_{l}^{K}(\Psi_{l}^{K})^{2}(I _{l}^{K})^{*}$ in the
following estimate.
\begin{multline*}\left\|P^{K}R^{K}\omega - R^{K}P\omega \right\|\leq
\sum_{l}\left\|I_{l}^{K}(\Psi_{l}^{K})^{2}(I
_{l}^{K})^{*}R^{K}\omega
-R^{K}i_{l}\psi_{l}^{2}i_{l}^{*}\omega\right\|+O(h)\\
\leq\sum_{l}\left(\left\|I_{l}^{K}(\Psi_{l}^{K})^{2}((I_{l}^{K})^{*}R^{K}\omega
-
R^{K}i_{l}^{*}\omega)\right\|+\left\|I_{l}^{K}((\Psi_{l}^{K})^{2}R^{K}i_{l}^{*}\omega
-R^{K}\psi_{l}^{2}i_{l}^{*}\omega)\right\|\right)\\
+\sum_{l}\left\|(R^{K}i_{l}-I_{l}^{K}R^{K})\psi_{l}^{2}i_{l}^{*}\omega\right\|+O(h)\leq
C_{i}h^{q-\frac{N}{2}+1}\|\omega\|_{r}^{2}.
\end{multline*}
Above  we used  (\ref{Iiest}) and (\ref{Istarest}) already proved
over each $\overline{U}_{l}$, and the inequality
\[\left\|(\Psi_{l}^{K})^{2}R^{K}\omega -
R^{K}\psi_{l}^{2}\omega\right\|^{2}\leq
C_{i}h^{2q-N+2}\|\omega\|_{r}^{2}\] which is proved in exactly the
same manner as (\ref{Iiest}). Thus (\ref{projR}) is shown to hold
for arbitrary $E$.

The estimate (\ref{projW}) can be easily derived from (\ref{projR}).
Finally, the equality (\ref{projRf}) follows directly from the
definitions.
\end{proof}

We use this lemma to deduce that the combinatorial connection
approximates the smooth one in the appropriate sense.

\begin{proposition}\label{mainest2}
For $r>\frac{N}{2}+1$, every  section $f \in \Omega^{0}(M,E)$ and
every 0-cochain $c \in C(K,E)$ one has

\begin{equation}\label{est1gen}\|\widetilde{W}^{K}\nabla ^{K}\widetilde{R}^{K}f-
\nabla f \|\leq C_{1}^{i}h\|f\|_{r}.\end{equation}

\end{proposition}

\begin{proof} We estimate as follows,
using (\ref{projRf}), (\ref{projW}) and (\ref{best}).
\begin{multline*}\|\widetilde{W}^{K}\nabla ^{K}\widetilde{R}^{K}f- \nabla f
\|=\|i^{*}W^{K}P^{K}d^{K}P^{K}R^{K}if-i^{*}dif\|\\
=\|i^{*}W^{K}P^{K}R^{K}dif-
i^{*}dif\|\leq\|(W^{K}P^{K}R^{K}-P^{K}W^{K}R^{K})dif\|\\
+\|W^{K}R^{K}dif-dif\|\leq C_{1}^{i}h\|f\|_{r}\end{multline*}
\end{proof}

The first inequality in Theorem \ref{convgeneral} follows as in the
proof of Theorem \ref{conver1}.

\subsection{Estimate from above}

In this subsection we estimate from above the eigenvalues of the
connection Laplacian in terms of the eigenvalues of its
discretizations, i.e. we prove the second inequality in Theorem
\ref{convgeneral}. This will follow from Proposition \ref{mainest2}
below. We will need the following strengthened version of Lemma
\ref{lebelow}.

\begin{lemma}\label{leabove}

For $r>\frac{N}{2}+2$ and every form $\omega \in \Omega(M,V)$  one
has
\begin{equation} \label{projWd}
\left\|dW^{K}P^{K}R^{K}\omega - dW^{K}R^{K}P\omega \right\|\leq
C_{i}h\|w\|_{r},
\end{equation}
\begin{equation} \label{projRd}
\left\|dW^{K}P^{K}R^{K}\omega - dPW^{K}R^{K}\omega \right\|\leq
C_{i}h\|w\|_{r}.
\end{equation}
Above $C_{i}$  is a  constant depending only on the map $i$ and its
derivatives.
\end{lemma}

\begin{proof} The argument is similar to the proof of Lemma
\ref{lebelow} but now  one has to use in addition the derivation
property of $d^{K}$. The details are omitted.
\end{proof}
We now use this lemma to derive the following proposition which is
the analog of part ii) of Proposition \ref{mainest1} in this
situation.
\begin{proposition}\label{mainest2}
For every 0-cochain $c \in C(K,E)$ one has
\begin{equation}\label{est2gen}\|\widetilde{W}^{K}\nabla ^{K}c - \nabla \widetilde{W}^{K}c
\|\leq C^{i}h (\|\widetilde{W}^{K}c\|+\|\nabla
\widetilde{W}^{K}c\|). \ee
\end{proposition}

\begin{proof}

We first prove the analogous estimate in the case when
$c=\widetilde{R}^{K}f$ for some smooth section $f$. Using
(\ref{projW}), (\ref{projRd}) and (\ref{projRf}) one has for $r$
large enough:
\begin{multline*}\|\widetilde{W}^{K}\nabla ^{K}\widetilde{R}^{K}f
 - \nabla \widetilde{W}^{K}\widetilde{R}^{K}f
\|=\|i^{*}W^{K}P^{K}R^{K}d\,if-i^{*}dPW^{K}P^{K}R^{K}if\| \leq \\
\|(W^{K}P^{K}R^{K}-PW^{K}R^{K})d\,if\|
+\|dW^{K}P^{K}R^{K}if-dPW^{K}R^{K}if\| \leq C_{2}^{i}h
\|f\|_{r}.\end{multline*} Now we note that (\ref{projRf}) implies
that $c=\widetilde{R}^{K}\widetilde{W}^{K}$ and  proceed exactly as
in the proof of part ii) of Proposition \ref{mainest1}. The only
difference now is that the image of $\widetilde{W}^{K}$
  does not consist of piecewise linear sections of $E$ anymore. However one can pass from
  the $r$-th to the first Sobolev norm on each $N$-simplex simply by a repeated use of the
  Leibniz rule.
\end{proof}

It remains to prove Lemma \ref{inject}. Observe that according to
(\ref{Iandi}) it suffices to prove that $i^{*}W^{K}i^{K}$ is
injective. \pagebreak Suppose that $c_{0}\in C^{q}(K,E)$ is a
nonzero cochain supported on a single $q$-simplex $\sigma$ with
reference point $p_{\sigma}$. A direct computation in local
coordinates  shows that $i^{*}W^{K}i^{K}c_{0}$ is not identically
zero. Now let $c$ be arbitrary nonzero cochain and suppose that
$i^{*}W^{K}i^{K}c=0$. Suppose also that $c$ is not zero at a simplex
$\sigma$. Then Proposition \ref{Wprop} will imply that
$j^{*}_{\sigma}(i^{*}W^{K}i^{K}c)=i^{*}j^{*}_{\sigma}(W^{K}i^{K}c)$
is not zero, a contradiction.

\section*{Acknowledgement}
This work is part of the author's Ph.D. thesis completed under the
guidance of Professor J\'ozef Dodziuk to whom many thanks are due
for his constant support.

\end{document}